\documentclass[12pt]{amsart}
\usepackage[english]{babel}
\usepackage{amsmath}
\usepackage{amsthm}
\usepackage{amssymb}
\usepackage{mathrsfs}
\usepackage{enumerate}
\usepackage[notcite, final, notref]{showkeys}
\usepackage{amsfonts}
\usepackage{dsfont}
\usepackage{tikz}
\usepackage[siunitx]{circuitikz}
\usetikzlibrary{arrows,calc,chains,shapes,dsp}
\def\C{\mathbb C}
\def\R{{\mathbb R}}
\newtheorem{Pa}{Paper}[section]
\newtheorem{Tm}[Pa]{{\bf Theorem}}
\newtheorem{La}[Pa]{{\bf Lemma}}
\newtheorem{Ob}[Pa]{{\bf Observation}}
\newtheorem{Cy}[Pa]{{\bf Corollary}}
\newtheorem{Rk}[Pa]{{\bf Remark}}
\newtheorem{Pn}[Pa]{{\bf Proposition}}

\newtheorem{Ex}[Pa]{{\bf Example}}
\newtheorem{Dn}[Pa]{{\bf Definition}}
\usepackage{xcolor}

\title[On the Hyper-Lyapunov Matrix Inclusion]
{On the Hyper-Lyapunov Matrix Inclusions}

\author[I. Lewkowicz]{Izchak Lewkowicz}
\address{(IL) School of Electrical and Computer Engineering\\
Ben-Gurion University of the Negev\\ P.O.B. 653\\ Beer-Sheva, 84105\\
Israel}
\email{izchak@bgu.ac.il}

\begin{document}
\bibliographystyle{plain}

\begin{abstract}
Gantmacher-Lyapunov Theorem (1950's) characterizes matrices whose
spectrum lies in the right-half of the complex plane. Here this
result is refined to Hyper-Lyapunov inclusion for matrices whose
spectrum lies in some disks within the right-half plane. These
disks turn to be closed under inversion, and when
their radius approaches infinity, the original result is recovered.
\vskip 0.2cm

\noindent
Hyper-Lyapunov inclusions are formulated through Quadratic Matrix
Inequalities and so are the analogous Hyper-Stein sets of matrices
whose spectrum lies within a {\em sub}-unit disk.
\vskip 0.2cm

\noindent
As a by-product, it is shown that these disks closed under inversion,
are a natural tool to understanding the Matrix Sign Function iteration
scheme, used in matrix computations.
\end{abstract}
\maketitle


\date{today}
\tableofcontents

\bibliographystyle{plain}
\section{Introduction}
\setcounter{equation}{0}

Gantmacher-Lyapunov Theorem (1950's) characterizes matrices whose
spectrum lies in the right-half of the complex plane, see.
\cite[Chapter 5]{Gant1971} and the Stein Theorem (1960's)
characterizes matrices whose spectrum lies in the unit disk, see
\cite{Stein1965} and \cite{Tau1964}. It is well known that these
theorems are inter-related through the Cayley transform, see
\cite{Ando2001}, \cite[Section 5.3]{LanRod1995}, \cite{Tau1964}.
\vskip 0.2cm

\noindent
By introducing a parameter $\eta\in(1,~\infty]$ we refine the
above to Hyper-Stein matrices: For a an arbitrary non-singular
Hermitian matrix $H$, we define the family
\[
\mathbf{S}_H({\scriptstyle\eta})=\left\{ A\in\C^{n\times n}~:~
\left({\scriptstyle\frac{\eta-1}{\eta+1}}H-A^*HA\right)\succ 0~\right\}.
\]
Whenever $H$ is positive definite, the spectrum of each of these
matrices lies in a disk centered at the origin, with radius of
${\scriptstyle\frac{\sqrt{\eta+1}}{\sqrt{\eta-1}}}$. This family
is closed under product among its elements. For details, see
Section \ref{sec:sub-unitDisk}.
\vskip 0.2cm

\noindent
Hyper-Lyapunov matrices are analogously defined:
For a an arbitrary non-singular
Hermitian matrix $H$, we define the family
\[
\mathbf{L}_H({\scriptstyle\eta})=\left\{ A\in\C^{n\times n}~:~
(HA+A^*H)\succ{\scriptstyle\frac{1}{\eta}}(A^*HA+H)~\right\}.
\]
Whenever $H$ is positive definite, the spectrum of each of these
matrices lies in a disk centered at $\eta+0i$, with radius of
$\sqrt{{\eta}^2-1}$ (and in particular within the right-half
plane), see Figure \ref{Figure:Eta=sqrt{2}}. This family
is closed under inversion.

\begin{figure}[ht!]
\begin{tikzpicture}[scale=2.2,cap=round]
    \tikzstyle{axes}=[]
    \tikzstyle{important line}=[very thick]
    \tikzstyle{information text}=[rounded corners,fill=red!10,inner sep=1ex]
 
   \draw[style=help lines,step=0.5cm];
     \begin{scope}[style=axes]
   \
       \draw[->] (-0.1,0) -- (2.7,0) node[right] {Real};
       \draw[->] (0,-1.1) -- (0,1.2) node[above] {Imaginary};
       \draw[-] (0.7071067,0.0)  -- (0.7071067,0.7071067) node[right]{};
       \draw[->] (0,0) -- (0.7071067,0.7071067) node[right]{};
    \foreach \x/\xtext in {0.7071067/{\scriptstyle\frac{1}{\eta}},1.41421356/{\scriptstyle\eta},
2.41421356/{\scriptstyle\eta+\sqrt{{\eta}^2-1}}}
 \draw[xshift=\x cm] (0pt,1pt) -- (0pt,-1pt) node[below,fill=white]
 {$\xtext$}; 
      \foreach \y/\ytext in {-1/{-{\scriptstyle\sqrt{{\eta}^2-1}}},
0.7071067/{\scriptstyle\frac{\sqrt{{\eta}^2-1}}{\eta}},
1/{\scriptstyle\sqrt{{\eta}^2-1}}}
        \draw[yshift=\y cm] (1pt,0pt) -- (-1pt,0pt) node[left,fill=white]
          {$\ytext$};
     \end{scope}
\draw[] (0,0) -- (3mm,0pt) arc(0:45:3mm);
\draw (22:2mm) node {$\scriptstyle\theta$};
   \draw[arrows=->,style=important line, red] (1.41421356,0) circle (1);
 \end{tikzpicture}
\caption{$\mathbb{D}_{\rm inv}(\eta),~~{\scriptstyle\sin(\theta)=
\frac{\sqrt{{\eta}^2-1}}{\eta}}$}
\label{Figure:Eta=sqrt{2}}
\end{figure}

\noindent
Here are some of the properties of this family.

\begin{Tm}\label{Tm:ClosureOfL_H(r)}
For prescribed matrix $~P\succ 0~$ and a scalar
\mbox{${\scriptstyle\eta}\in(1,~\infty]$}, the following is true.
\begin{itemize}

\item[(i)~~~]{}
The open set $~\mathbf{L}_P({\scriptstyle\eta})$ in Eq.
\eqref{eq:DefLh(eta)} is convex and closed under inversion.
\vskip 0.2cm

\item[(ii)~~]{}
For ${\scriptstyle\eta}=\infty$ the set $~\mathbf{L}_P$ is a
maximal convex set of matrices whose spectrum is in $\C_R~$.
\vskip 0.2cm

\item[(iii)~]{}
For $~P=I$, the set $~{\mathbf L}_I({\scriptstyle\eta})$, is 
in addition matrix-convex.
\vskip 0.2cm

\item[(iv)~]{} For arbitrary non-singular matrix $T$,
\begin{equation}\label{eq:Similarity}
T^{-1}\mathbf{L}_P({\scriptstyle\eta})T=
\mathbf{L}_{T^*PT}({\scriptstyle\eta}).
\end{equation}
\end{itemize}
\end{Tm}

\noindent
For further details, see Subsection \ref{Subsect:Hyper-Lyapunov}.
\vskip 0.2cm

\noindent
Expectedly, the families $\mathbf{S}_H({\scriptstyle\eta})$ and
$\mathbf{L}_H({\scriptstyle\eta})$ are inter-related through the Cayley
transform. Moreover, when $\eta~\rightarrow~\infty$, the classical
results are recovered.
\vskip 0.2cm

\noindent
Following the original results, the Hyper-Lyapunov and the
Hyper-Stein criterion are applied to continuous-time and
discrete-time stability problems, respectively. Consequently,
the smaller $\eta$ is, the estimate on the time-wise
trajectory improves.
\vskip 0.2cm

\noindent
In passing, in Subsection \ref{Sect:MSF} it is pointed out that
these invertible disks are natural (and visual) objects in
demonstrating the convergence of the Matrix Sign Function
iteration scheme, used in matrix computations, see e.g.
\cite[Chapter 5]{Higham2008}, \cite{KennLaub1995},
\cite[Chapter 22]{LanRod1995} and also \cite{CohenLew1997a},
\cite{CohenLew2007}, \cite{Lewk1999a}, \cite{LewRodYar2005}.

\begin{figure}[ht!]
 \begin{tikzpicture}[scale=0.75,cap=round]
    \tikzstyle{axes}=[]
    \tikzstyle{important line}=[very thick]
    \tikzstyle{information text}=[rounded corners,fill=red!10,inner sep=1ex]
    \begin{scope}[style=axes]
  \
      \draw[->] (-7.5,0) -- (8.3,0) node[right] {\rm Real};
      \draw[->] (0,-4.1) -- (0,4.3) node[above] {\rm Imaginary};

      \foreach \x/\xtext in {
-7.3276/{\scriptstyle-2-\sqrt{3}-\sqrt{6+4\sqrt{3}}},
-3.7321/{\scriptstyle-2-\sqrt{3}},
-2/{\scriptstyle-2},
1.25/{\scriptstyle\frac{5}{4}},
3.7321/{\scriptstyle 2+\sqrt{3}},
5/{\scriptstyle 5},
6/{\scriptstyle 6},
7/{\scriptstyle 7},
8/{\scriptstyle 8}
}
\draw[xshift=\x cm] (0pt,1pt) -- (0pt,-1pt) node[below,fill=white]
              {$\xtext$}; 
      \foreach \y/\ytext in {
-3.5956/{\scriptstyle-\sqrt{6+4\sqrt{3}}},
-1.7321/-{\scriptstyle\sqrt{3}},
2/{\scriptstyle 2},
3/{\scriptstyle 3},
4/{\scriptstyle 4}
}
       \draw[yshift=\y cm] (1pt,0pt) -- (-1pt,0pt) node[left,fill=white]
              {$\ytext$};
    \end{scope}
   \draw[arrows=->,style=important line, blue] (-3.7321,0) circle (3.6);
   \draw[arrows=->,style=important line, blue] (3.7321,0) circle (3.6);
   \draw[arrows=->,style=important line, red] (-2,0) circle (1.7321);
   \draw[arrows=->,style=important line, red] (2,0) circle (1.7321);
   \draw[arrows=->,style=important line, green] (-1.25,0) circle (0.75);
   \draw[arrows=->,style=important line, green] (1.25,0) circle (0.75);
   \draw[arrows=->,style=important line, gray] (-1.025,0) circle (0.225);
   \draw[arrows=->,style=important line, gray] (1.025,0) circle (0.225);
 \end{tikzpicture}
\caption{
\mbox{$\mathbb{D}_{\rm inv}({\scriptstyle\eta}):~
{\scriptstyle\eta}={\scriptstyle 2+\sqrt{3}}~-{\rm blue},~~~
{\scriptstyle h_1(\eta)}=2~-{\rm red},~~~
{\scriptstyle h_2(\eta)}={\scriptstyle\frac{5}{4}}~-{\rm green},~~
{\scriptstyle h_3(\eta)}={\scriptstyle\frac{41}{40}}~-{\rm gray}$}}
\label{Figure:MSF}
\end{figure}

\noindent
We end this work by showing, in Subsection
\ref{SubSec:UnifyingFramework}, that the (Hyper)-Stein and the
(Hyper)-Lyapunov inclusions can all be casted in a the same
Quadratic Matrix Inequality framework.

\section{Background}
\setcounter{equation}{0}
\label{sec:Bacground}

In this section we summarize the background relevant to the sequel.
\vskip 0.2cm

\noindent
Let $\mathbf{H}_n$ be the set of non-singular $n\times n$ Hermitian
matrices. For an Hermitian matrix $M$, we shall use the convention that 
\mbox{$M\succ 0$} \mbox{($M\succcurlyeq 0$)} means that $M$ is positive
(semi)-definite.
\vskip 0.2cm

\noindent
For a prescribed $H\in\mathbf{H}_n$ let us define the set of all
matrices satisfying the following Lyapunov inclusion,
\begin{equation}\label{eq:OriginalLyapunov}
\begin{matrix}
\mathbf{L}_H=&\left\{ A\in\C^{n\times n}~:~HA+A^*H\succ 0~\right\}
\\~\\
\overline{\mathbf L}_H=&\left\{ A\in\C^{n\times n}~:~
HA+A^*H\succcurlyeq 0~\right\}.
\end{matrix}
\end{equation}
We consider $\overline{\mathbf L}_H$ as the closure of the open set
$\mathbf{L}_H$, see e.g. \cite{CohenLew1997a}.
\vskip 0.2cm

\noindent
We also have an analogous family of the Stein inclusion,
\begin{equation}\label{eq:OriginalStein}
\begin{matrix}
\mathbf{S}_H=&\left\{ A\in\C^{n\times n}~:~
H-A^*HA\succ 0~\right\}
\\~\\
\overline{\mathbf S}_H=&\left\{ A\in\C^{n\times n}~:~H-A^*HA\succcurlyeq 0
~\right\}.
\end{matrix}
\end{equation}
As before, $\overline{\mathbf S}_H$ is the closure of the open set
$\mathbf{S}_H$. A characterization of the set $\mathbf{S}_H$, through
its structure was introduced by T. Ando in \cite[Theorem 3.5]{Ando2004}.
\vskip 0.2cm

\noindent
We next resort to the Cayley transform.

\begin{Dn}\label{Dn:MatrixCayleyTransform}
{\rm
We denote by $\mathcal{C}(A)$ the Cayley transform of a matrix
$A\in\C^{n\times n}$, \mbox{$-1\not\in{\rm spec}(A)$},
\begin{center}
$\begin{matrix}
\mathcal{C}\left(A\right)&:=&\left(I_n-A\right)\left(I_n+A\right)^{-1}
=-I_n+2\left(I_n+A\right)^{-1}
\end{matrix}$
\end{center}
\qed
}
\end{Dn}

\noindent
Recall that the Cayley transform is involutive in the sense that,
whenever well defined,
\[
\mathcal{C}\left(\mathcal{C}\left(A\right)\right)=A.
\]
The following is classical, see e.g. \cite{Ando2004}, 
\cite{Stein1965}, \cite{Tau1964}.

\begin{Pn}\label{Pn:CayleyLyapStein}
For a given $H\in\mathbf{H}_n$ let $\mathbf{L}_H$ and $\mathbf{S}_H$
(and the respective closures) be as in Eqs. \eqref{eq:OriginalLyapunov}
and \eqref{eq:OriginalStein} respectively. Then one has that
\begin{equation}\label{eq:CayleyStein}
\mathcal{C}\left(\mathbf{S}_H\right)=\mathbf{L}_H~.
\end{equation}
\end{Pn}
\vskip 0.2cm

\noindent
We next examine the structure of the two sets in Eqs.
\eqref{eq:OriginalLyapunov}, \eqref{eq:OriginalStein} and
\eqref{eq:CayleyStein}. To this end, we need the following.

\begin{Dn}\label{Dn:MatrixConvex}
{\rm
A family $\mathbf{A}$, of square matrices (of various dimensions) is
said to be~ {\em matrix-convex}, if for all natural $k$, $n$,
\begin{equation}\label{eq:Isometry1}
\sum\limits_{j=1}^k\upsilon_j^*\upsilon_j=I_n
\quad\quad
\begin{smallmatrix}
\forall{\upsilon}_j\in\C^{{\eta}_j\times n}
\\~\\
\forall\eta_j~,
\end{smallmatrix}
\end{equation}
one has that having $A_1,~\ldots~,~A_k$ (of dimensions
\mbox{${\eta}_1\times{\eta}_1$} through
\mbox{${\eta}_k\times{\eta}_k$}) within
$\mathbf{A}$, implies that also the $n\times n$ matrix
\[
\sum\limits_{j=1}^k\upsilon_j^*A_j\upsilon_j
\]
belongs to $\mathbf{A}$.
\qed
}
\end{Dn}
\vskip 0.2cm

\noindent
For background on matrix-convexity, see e.g. \cite{EffrWink1997} and
more recently, \cite{EverHeltKlepMcCull2018}, \cite{Kriel2019} and
\cite{PassShalSol2018}. See also \cite[Sections 2,3]{Lewk2020a}.
\vskip 0.2cm

\noindent
The set $\mathbf{L}_H$ where $H\in\mathbf{H}_n$, was first explored 
in \cite{CohenLew1997a}. Here we summarize its main properties.

\begin{Tm}\label{Tm:L_I=P+iH}\cite[Theorems 2.4 and 4.1]{Lewk2020a}
The following statements are true.

\begin{itemize}
\item[(i)~~~]{}
The set $\mathbf{L}_H$ where $H\in\mathbf{H}_n$, is
a cone closed under inversion. It is maximal open convex
set of non-singular matrices, containing the matrix $H$.
\vskip 0.2cm

\item[(ii)~~]{}
The set $\mathbf{L}_I$ is in addition matrix-convex.

\item[(iii)~]{}
Conversely, a cone closed under inversion and a maximal open
matrix-convex set of non-singular matrices, containing the
matrix $I$, is the set ${\mathbf L}_I$.
\end{itemize}
\end{Tm}

\noindent
Recall that items (ii) and (iii) refer to all matrices $A$ 
where $A+A^*\succ 0$, of various dimensions.
\vskip 0.2cm

\noindent
Similarly, here are the main properties of $\mathbf{S}_H$.

\begin{Tm}\label{Tm:Stein} \cite[Theorem 2.1]{Lewk2020c}
For $H\in\mathbf{H}_n$ the set $\mathbf{S}_H$ is open, convex,
closed under multiplication by $c\in\C$, where $1\geq |c|$
and under product among its elements, i.e. whenever 
$A\in{\scriptstyle\frac{1}{\alpha}}\mathbf{S}_H$ and
$B\in{\scriptstyle\frac{1}{\beta}}\mathbf{S}_H$, for some
$\alpha$, $\beta>0$, then their product satisfies
$AB\in{\scriptstyle\frac{1}{\alpha\beta}}\mathbf{S}_H~$. 
\vskip 0.2cm

\noindent
When $H=I$, the set $\mathbf{S}_I$ is matrix-convex,
\end{Tm}

\noindent
To gain intuition, we temporarily confine the discussion to
scalar set-up.

\section{Disks}\label{Sec:Disks}
\setcounter{equation}{0}
\label{sec:Disks}

\subsection{Sub-Unit ~Disks}
\label{Sect:SubUnitDisk}

\noindent
We shall use the following notation for disks in the complex plane 
\[
\begin{matrix}
{\mathbb D}({\scriptstyle{\rm Center}},~{\scriptstyle{\rm Radius}})
=\{ {\scriptstyle x}\in\C~:~{\scriptstyle{\rm Radius}}>
|{\scriptstyle x}-{\scriptstyle{\rm Center}}|~\}
\\~\\
\overline{\mathbb D}({\scriptstyle{\rm Center}},~{\scriptstyle{\rm Radius}})
=\{ {\scriptstyle x}\in\C~:~{\scriptstyle{\rm Radius}}\geq
|{\scriptstyle x}-{\scriptstyle{\rm Center}}|~\}
\end{matrix}
\quad\quad\quad\quad\begin{matrix}
{\rm Center}\in\C\\~\\ ~{\rm Radius}>0.
\end{matrix}
\]

\noindent
A key-parameter in this work is $\eta$, which assumes values in
$(1,~\infty]$. We start describing sub-unit disks of the form,
\begin{equation}\label{eq:Dorigin}
\mathbb{D}_{\rm origin}\left(\eta\right):=
\mathbb{D}(\underbrace{0}_{\rm Center},~
\underbrace{\scriptstyle\sqrt{\frac{\eta-1}{\eta+1}}}_{\rm Radius})
\quad\quad\quad\eta\in(1,~\infty].
\end{equation}
This is illustrated in Figure \ref{Figure:SubUnitDisks}.\quad
We next examine properties of these disks.

\begin{figure}[ht!]
 \begin{tikzpicture}[scale=6.0,cap=round]
    \tikzstyle{axes}=[]
    \tikzstyle{important line}=[very thick]
    \tikzstyle{information text}=[rounded corners,fill=red!10,inner sep=1ex]
    \begin{scope}[style=axes]
  \
      \draw[->] (-0.6,0) -- (0.65,0) node[right] {Real};
      \draw[->] (0,-0.6) -- (0,0.6) node[above] {Imaginary};

      \foreach \x/\xtext in 
{- 0.5/-{\scriptstyle\frac{\sqrt{\eta-1}}{\sqrt{\eta+1}}},
-0.25/-{\scriptstyle\frac{\eta-1}{\eta+1}},
0.25/{\scriptstyle\frac{\eta-1}{\eta+1}},
0.5/{\scriptstyle\frac{\sqrt{\eta-1}}{\sqrt{\eta+1}}}
}
\draw[xshift=\x cm] (0pt,1pt) -- (0pt,-1pt) node[below,fill=white]
              {$\xtext$}; 

      \foreach \y/\ytext in 
{- 0.5/-{\scriptstyle\frac{\sqrt{\eta-1}}{\sqrt{\eta+1}}},
-0.25/-{\scriptstyle\frac{\eta-1}{\eta+1}},
0.25/{\scriptstyle\frac{\eta-1}{\eta+1}},
0.5/{\scriptstyle\frac{\sqrt{\eta-1}}{\sqrt{\eta+1}}}
}
       \draw[yshift=\y cm] (1pt,0pt) -- (-1pt,0pt) node[left,fill=white]
              {$\ytext$};
    \end{scope}
   \draw[arrows=->,style=important line, red] (0,0) circle (0.5);
   \draw[arrows=->,style=important line, blue] (0,0) circle (0.25);
 \end{tikzpicture}
\caption{$
\mathbb{D}_{\rm origin}\left({\scriptstyle\frac{\sqrt{\eta-1}}
{\sqrt{\eta+1}}}\right)~-{\scriptstyle{\rm red}},\quad
\mathbb{D}_{\rm origin}\left({\scriptstyle\frac{\eta-1}{\eta+1}}
\right)~-{\scriptstyle{\rm blue}},\quad\eta\in(1,~\infty].$}
\label{Figure:SubUnitDisks}
\end{figure}

\begin{La}\label{La:ProductHBeta}
For arbitrary scalars
$~{\scriptstyle\eta}_a ,~{\scriptstyle\eta}_b\in(1,~\infty]~$ let
\begin{equation}\label{eq:Eta1Eta2}
{\eta}_c=\frac{{\scriptstyle 1}+{\scriptstyle\eta}_a{\scriptstyle\eta}_b}
{{\scriptstyle\eta}_a+{\scriptstyle\eta}_b}~.
\end{equation}
Then the following is true.
\begin{itemize}

\item[(i)~~~]{}
Let ${\eta}_c$ be as in Eq. \eqref{eq:Eta1Eta2}, then
\[
{\eta}_c={\scriptstyle\theta}{\eta}_a+
(1-{\scriptstyle\theta}){\scriptstyle\frac{1}{{\eta}_a}}\quad\quad\quad
{\scriptstyle\theta}={\scriptstyle\frac{{\eta}_b}{{\eta}_a+{\eta}_b}}~~.
\]
Alternatively,
\[
\frac{1}{{\eta}_c}=\frac{1}{{\eta}_a+{\frac{1}{{\eta}_b}}}
+
\frac{1}{{\eta}_b+{\frac{1}{{\eta}_a}}}~.
\]
\item[(ii)~~]{} For
$~{\eta}_b\geq{\eta}_a~$, let 
$~{\eta}_c$ be as in Eq. \eqref{eq:Eta1Eta2}. Then,
\[
{\scriptstyle\frac{1}{2}}\left({\eta}_b+{\scriptstyle\frac{1}{{\eta}_b}}\right)
\geq{\eta}_c\geq
{\scriptstyle\frac{1}{2}}\left({\eta}_a+{\scriptstyle\frac{1}{{\eta}_a}}\right).
\]
\vskip 0.2cm

\item[(iii)~]{} A product of a pair of disks is given by,
\[
\mathbb{D}_{\rm origin}\left({\scriptstyle\frac{\sqrt{{\eta}_a-1}}
{\sqrt{{\eta}_a+1}}}\right)\cdot
\mathbb{D}_{\rm origin}\left({\scriptstyle\frac{\sqrt{{\eta}_b-1}}
{\sqrt{{\eta}_b+1}}}\right)=
\mathbb{D}_{\rm origin}\left({\scriptstyle\frac{\sqrt{{\eta}_c-1}}
{\sqrt{{\eta}_c+1}}}\right),
\]
with $~\eta_c$ as in Eq. \eqref{eq:Eta1Eta2}.
\end{itemize}
\end{La}

\noindent
{\bf Proof :}\quad (i)~ This is immediate.
\vskip 0.2cm

\noindent
(ii)~ Bounding ${\eta}_c$ from above and from below,
\[
\begin{matrix}
{\eta}_b({\eta}_b^2-1)&\geq&{\eta}_a({\eta}_b^2-1)&&~{\rm by~assumption}\\~\\
{\eta}_a+{\eta}_b^3&\geq&{\eta}_a{\eta}_b^2+{\eta}_b&&{\rm re-arranging~terms}\\~\\
({\eta}_b^2+1)({\eta}_a+{\eta}_b)&\geq&2{\eta}_b({\eta}_a{\eta}_b+1)
&&{\rm adding~the~same~terms~to~both~sides}\\~\\
{\scriptstyle\frac{1}{2}}({\eta}_b+\frac{1}{{\eta}_b})&\geq&{\eta}_c&&
{\rm dividing~each~side~by~}2{\eta}_b({\eta}_a+{\eta}_b)
\end{matrix}
\]
so the right-hand side is obtained. Similarly, for the
left-hand side,
\[
\begin{matrix}
{\eta}_b({\eta}_a^2-1)&\geq&{\eta}_a({\eta}_a^2-1)&&~{\rm by~assumption}\\~\\
{\eta}_a^2{\eta}_b+{\eta}_a&\geq&{\eta}_a^3+{\eta}_b&&{\rm re-arranging~terms}\\~\\
2{\eta}_a({\eta}_a{\eta}_b+1)&\geq&({\eta}_a^2+1)({\eta}_a+{\eta}_b)&&
{\rm adding~the~same~terms~to~both~sides}\\~\\
{\eta}_c&\geq&{\scriptstyle\frac{1}{2}}({\eta}_a+\frac{1}{{\eta}_a})&&
{\rm dividing~each~side~by~}2{\eta}_a({\eta}_a+{\eta}_b),
\end{matrix}
\]
establishes, this part of the claim.
\vskip 0.2cm

\noindent
(iii)~ In principle, for all positive ${\rm Radius}_a~, {\rm Radius}_b$~,
a product of disks centered at the origin satisfies,
\[
\mathbb{D}_{\rm origin}({\scriptstyle{\rm Radius}_{a})}\cdot\mathbb{D}_{\rm origin}(
{\scriptstyle{\rm Radius}_b})=\mathbb{D}_{\rm origin}({\scriptstyle{\rm Radius}_c}
),\quad{\rm with}\quad{\scriptstyle {\rm Rradius}_c}=
{\scriptstyle{\rm Radius}_a}\cdot{\scriptstyle {\rm Radius}_b}~,
\]
and in particular,
\[
\underbrace{\scriptstyle\frac{\sqrt{{\eta}_a-1}}{\sqrt{{\eta}_a+1}}}_{{{\rm Radius}}_a}
\underbrace{\scriptstyle\frac{\sqrt{{\eta}_b-1}}{\sqrt{{\eta}_b+1}}}_{{{\rm Radius}}_b}
=
\underbrace{\scriptstyle\frac{\sqrt{{\eta}_c-1}}{\sqrt{{\eta}_c+1}}}_{{{\rm Radius}}_c},
\]
which (after taking squares) implies that $\underbrace{\scriptstyle\left(1-
\frac{2}{{\eta}_a+1}\right)}_{{{\rm Radius}}_a^2}\underbrace{\scriptstyle
\left(1-\frac{2}{{\eta}_b+1}\right)}_{{{\rm Radius}}_b^2}
=
\underbrace{\scriptstyle\left(1-\frac{2}{{\eta}_c+1}\right)}_{{{\rm Radius}}_c^2}~$,
and in turn,
\[
{\scriptstyle\frac{{\eta}_a+{\eta}_b}{({\eta}_a+1)({\eta}_b+1)}}
={\scriptstyle\frac{1}{{\eta}_c+1}}~,
\]
so Eq. \eqref{eq:Eta1Eta2} is obtained.
\qed
\vskip 0.2cm

\noindent
Item (i) says that ${\eta}_c$ may be formulated in terms of
${\eta}_a$ and ${\eta}_b$ in two equivalent ways:
(a) a convex combination or (b) a special kind of a double
harmonic mean.
\vskip 0.2cm

\noindent
Recall that in item (iii) of Lemma \ref{La:ProductHBeta} we described
the product of a pair of sub-unit disks. Of a special interest is the
case where their radii are equal (${\eta}_a={\eta}_b={\eta}$), i.e.
when Eq. \eqref{eq:Eta1Eta2} takes the form,
\begin{equation}\label{eq:Eta1}
{{\eta}_c}_{|_{{\eta}_a={\eta}_b=\eta}}=\underbrace{
{\scriptstyle\frac{1}{2}}\left(\eta+{\scriptstyle\frac{1}{\eta}}
\right).}_{h_1(\eta)~{\rm see~Eq.}~\eqref{eq:h_j}}
\end{equation}
For example, in Figure \ref{Figure:SubUnitDisks} the blue disk
is a result of the product of a pair of identical red disks.
\vskip 0.2cm

\noindent
In Section \ref{sec:sub-unitDisk} below, we extend the above discussion
to matricial framework:
\[
\begin{matrix}
{\mathbb D}_{\rm origin}(\eta)~~{\rm Eq.~\eqref{eq:Dorigin}}
&\longrightarrow &\left\{A\in\C^{n\times n}~:~\sqrt{\frac{\eta-1}
{\eta+1}}>\| A\|_2~\right\}
\end{matrix}
\]

\subsection{Disks closed under inversion}
\label{Sect:InvertibleDisk}

In this subsection we exploit the previous disks to describe special
disks in the open left and right half planes, $\C_L$ and $\C_R$. To
this end we find it convenient to introduce the following disks,
where the subscript abbreviates ~{\em invertible} (meaning ``closed
under inversion"), see Figures \ref{Figure:BasicInvertibleDisks}
and \ref{Figure:MSF}.
\begin{equation}\label{eq:InvertibleDisk}
\mathbb{D}_{\rm inv}\left(\eta\right):=
\mathbb{D}(\underbrace{\scriptstyle i\cdot 0+\eta}_{\rm Center},~
\underbrace{\scriptstyle\sqrt{\eta^2-1}}_{\rm Radius})
\quad\quad\quad\eta\in(1,~\infty],
\end{equation}
As always, \mbox{$-\mathbb{D}_{\rm inv}\left(\eta\right)=
\mathbb{D}(\underbrace{\scriptstyle i\cdot 0-\eta}_{\rm Center},~
\underbrace{\scriptstyle\sqrt{\eta^2-1}}_{\rm Radius})\subset\C_L$},
thus without loss of generality, in the sequel we can focus on
$~\mathbb{D}_{\rm inv}$, disks in $\C_R~$.\quad Here are some
properties of these disks.

\begin{La}\label{La:Dinv}
For a parameter ${\scriptstyle\eta}\in(1,~\infty]$, let 
$~\mathbb{D}_{\rm inv}({\scriptstyle\eta})$ be as in Eq.
\eqref{eq:InvertibleDisk}.\\
Then the following statements are true.

\begin{itemize}

\item[(i)~~~]{}$\mathbb{D}_{\rm inv}\left(\eta\right)=\left\{
{\scriptstyle c}\in\C_R~:~{\scriptstyle\eta}>{\scriptstyle
\frac{|c|^2+1}{2{\rm Re}(c)}}~\right\}$.
\vskip 0.2cm

\item[(ii)~~]{}Through inversion, the disk $~\mathbb{D}_{\rm inv}
\left(\eta\right)$ is mapped onto itself, i.e. \mbox{$\left(
\mathbb{D}_{\rm inv}\left(\eta\right)\right)^{-1}=
\mathbb{D}_{\rm inv}\left(\eta\right)$.}
\vskip 0.2cm

\item[(iii)~]{}\mbox{$\mathcal{C}\left(\mathbb{D}_{\rm inv}
\left(\eta\right)\right)=\mathbb{D}_{\rm origin}(\eta)$}.
\end{itemize}
\end{La}

\noindent
Indeed, item (i) follows from Eq. \eqref{eq:InvertibleDisk} where
for some ${\scriptstyle\eta}\in(1,~\infty]$,
\[
\begin{matrix}
\mathbb{D}_{\rm inv}({\scriptstyle\eta})&=\{~{\scriptstyle c}\in\C_R~
:&{\scriptstyle\eta}^2-1>|{\scriptstyle c}-{\scriptstyle\eta}|^2&\}
\\~\\~&=\{
{\scriptstyle c}\in\C_R~:&{\scriptstyle\eta}^2-1>{\scriptstyle\eta}^2-
2{\scriptstyle\eta}{\scriptstyle\rm Re}({\scriptstyle c})
+|{\scriptstyle c}|^2&\}
\\~\\~&=\{
{\scriptstyle c}\in\C_R~:&2{\scriptstyle\eta}{\scriptstyle\rm Re}
({\scriptstyle c})>{\scriptstyle\eta}^2+1&\}
\\~\\~&=\{
{\scriptstyle c}\in\C_R~:&{\scriptstyle\eta}>{\scriptstyle
\frac{|c|^2+1}{2{\rm Re}(c)}}&\}.
\end{matrix}
\]
Item (ii) is obtained by substituting \mbox{${\scriptstyle c}~
\longrightarrow~{\scriptstyle\frac{1}{c}}$}, in item (i).
\vskip 0.2cm

\noindent
As to item (iii), we exploit the fact that this map is conformal
and first show that a boundary of the closure is mapped to a
boundary of a closure, i.e.
\[
\begin{matrix}
\mathcal{C}(\partial\overline{\mathbb D}_{\rm origin}
(\sqrt{\scriptstyle\frac{\eta-1}{\eta+1}}))
&=
\{\mathcal{C}
\left(\sqrt{\scriptstyle\frac{\eta-1}{\eta+1}}\cdot{e^{i\theta}}
\right)~:&{\scriptstyle\theta}\in[{\scriptstyle 0,~2\pi})\}
\\~\\~&=\{{\scriptstyle c}=\frac
{1-\sqrt{\scriptstyle\frac{\eta-1}{\eta+1}}\cdot{e^{i\theta}}}
{1+\sqrt{\scriptstyle\frac{\eta-1}{\eta+1}}\cdot{e^{i\theta}}}
~:&{\scriptstyle\theta}\in[{\scriptstyle 0,~2\pi})\}
\\~\\~&=\{{\scriptstyle c}=\frac{1-i\sqrt{\scriptstyle{\eta}^2-1}
\sin(\theta)}
{{\scriptstyle\eta}+\left(\sqrt{\scriptstyle{\eta}^2-1}\right)\cos(\theta)}
~:&{\scriptstyle\theta}\in[{\scriptstyle 0,~2\pi})\}
\end{matrix}
\]
and thus $~{\scriptstyle{\rm Re}(c)}=\frac{1}
{{\scriptstyle\eta}+\left(\sqrt{\scriptstyle{\eta}^2-1}\right)\cos(\theta)}~$
and $~|{\scriptstyle c}|^2+1=2{\scriptstyle\eta{\rm Re}(c)}$,
so indeed $~\frac{|c|^2+1}{2{\rm Re}(c)}=\eta$, which is on
$\partial\overline{\mathbb D}_{\rm inv}(\eta)$.
\quad
All is left is to show that interior is mapped to interior. Indeed
as $\mathcal{C}(0)=1+0i$, where $0\in\mathbb{D}_{\rm origin}$ and
$1+0i\in\mathbb{D}_{\rm inv}$, so the proof is complete.
\qed
\vskip 0.2cm

\noindent
As a side remark, recall that, Andres Rantzer already addressed (in the
context of weak Kharitonov Theorem) convex sets closed under inversion, 
see \cite{Rant1993} and for perspective \cite{Barmish1994}.
\vskip 0.2cm

\noindent
In Corollary \ref{Cy:QuantitativeHyperPositive1} below, we extend
the above discussion to matricial framework:
\[
\begin{matrix}
{\mathbb D}_{\rm inv}(\eta)~~{\rm Lemma~\ref{La:Dinv},~ item~(i)}
&\longrightarrow&\{A\in\mathbf{L}_I~:~\eta>\rho\left(
{\scriptstyle(A^*A+I)(A+A^*)^{-1}}\right)~\}
\end{matrix}
\]
Here is our motivation to introducing these disks in the context
of the Newton iterations (a.k.a. ``half iterations").

\begin{Pn}\label{Pn:HalfIterationDinv}
For $\eta\in(1,~\infty]$, the following is true
\[
c\in\mathbb{D}_{\rm inv}(\eta)\quad\Longleftrightarrow\quad
{\scriptstyle\frac{1}{2}}\left(c+{\scriptstyle\frac{1}{c}}
\right)\in\mathbb{D}_{\rm inv}(
{\scriptstyle\frac{1}{2}}(\eta+
{\scriptstyle\frac{1}{\eta}})).
\]
In $\C_L$, an analogous statement holds for
$~-\mathbb{D}_{\rm inv}(\eta)$.
\end{Pn}

\noindent
{\bf Proof :}\quad We show it in two ways. Indirect: Combine
item (iii) in Lemma \ref{La:ProductHBeta}
(see Eq. \eqref{eq:Eta1}) together with item (iii) of Lemma \ref{La:Dinv}
(see also item(ii) of Lemma \ref{La:CayleySignInversion}).
below.
\vskip 0.2cm

\noindent
Direct: We again exploit the fact that boundary is mapped to boundary.
Recall,
\[
\partial\overline{\mathbb D}_{\rm inv}(\eta)\left\{ c\in\C_r~:~
{\scriptstyle\frac{|c|^2+1}{2{\rm Re}(c)}}=\eta~\right\},
\]
denote $c_1:=
{\scriptstyle\frac{1}{2}\left(c+\frac{1}{c}\right)}$,
and examine the following quantity
\[
\begin{matrix}
\frac{|c_1|^2+1}{2{\rm Re}(c_1)}
&=&
\frac{\left({\rm Re}(c_1)\right)^2+\left({\rm Im}(c_1)\right)^2+1}{2{\rm Re}(c_1)}
=
\frac{
{\scriptstyle\left(
\frac{1}{2}{\rm Re}(c)\right)^2\left(1+\frac{1}{|c|^2}\right)^2}
+
{\scriptstyle\left(\frac{1}{2}{\rm Im}(c)\right)^2\left(1-\frac{1}{|c|^2}\right)^2}
+1}
{2
{\scriptstyle\frac{1}{2}{\rm Re}(c)\left(1+\frac{1}{|c|^2}\right)}}
\\~\\~&=&\frac{{\scriptstyle\left(\frac{1}{2}{\rm Re}(c)\right)^2}
{\scriptstyle\left(\left(1-\frac{1}{|c|^2}\right)^2+
\frac{4}{|c|^2}
\right)}
+
{\scriptstyle\left(\frac{1}{2}{\rm Im}(c)\right)^2\left(1-\frac{1}{|c|^2}\right)^2}
+1}
{2
{\scriptstyle\frac{1}{2}{\rm Re}(c)\left(1+\frac{1}{|c|^2}\right)}
}
=
\frac
{\scriptstyle\left(\frac{|c|}{2}\right)^2\left(1-\frac{1}{|c|^2}\right)^2+
\left(\frac{1}{|c|}{\rm Re}(c)\right)^2+
1}
{\scriptstyle{\rm Re}(c)\left(1+\frac{1}{|c|^2}\right)}
\\~\\~&=&
\frac{\scriptstyle\left(\frac{|c|^2}{2}\right)^2\left(1-\frac{1}{|c|^2}
\right)^2+\left({\rm Re}(c)\right)^2+|c|^2}
{\scriptstyle{\rm Re}(c)\left(|c|^2+1\right)}
=
\frac{\scriptstyle\frac{1}{4}\left(|c|^2+1\right)^2+\left({\rm Re}(c)\right)^2}
{\scriptstyle{\rm Re}(c)\left(1+|c|^2\right)}={\scriptstyle\frac{1}{2}}
\left(\frac{\scriptstyle|c|^2+1}{\scriptstyle 2{\rm Re}(c)}
+\frac{\scriptstyle2{\rm Re}(c)}{\scriptstyle 1+|c|^2}\right)
\\~\\~&=&
{\scriptstyle\frac{1}{2}}
\left(\eta+\scriptstyle\frac{1}{\eta}\right).
\end{matrix}
\]
Namely,
$c_1\in\partial\overline{\mathbb D}_{\rm inv}\left({\scriptstyle\frac{1}{2}(\eta
+\frac{1}{\eta})}\right)$, so indeed boundary is mapped to boundary.
Finally as $1+0i$ is an invariant interior point of both disks, the
construction is complete.
\qed
\vskip 0.2cm

\noindent
In Observation \ref{Ob:DisksIterations} below we further exploit
Observation \ref{Pn:HalfIterationDinv}.
\vskip 0.2cm

\begin{Ex}\label{Ex:DinvSquare}
{\rm 
In Figure \ref{Figure:SubUnitDisks}, the blue disk is the (point-wise) square
of the red disk. 
\vskip 0.2cm

\noindent
In Figure \ref{Figure:BasicInvertibleDisks} blue disk is
the (point-wise) image of the red disk through the map ${\scriptstyle
\frac{1}{2}(x+\frac{1}{x})}$.
\vskip 0.2cm

\noindent
The red (blue) disk in Figure \ref{Figure:BasicInvertibleDisks} is the image,
through the Cayley transform, of the red (blue) disk from Figure
\ref{Figure:SubUnitDisks}.
}
\qed
\end{Ex}

\begin{figure}[ht!]
\begin{tikzpicture}[scale=2.0,cap=round]
    \tikzstyle{axes}=[]
    \tikzstyle{important line}=[very thick]
    \tikzstyle{information text}=[rounded corners,fill=red!10,inner sep=1ex]
 
   \draw[style=help lines,step=0.5cm];
     \begin{scope}[style=axes]
   \
       \draw[->] (-0.3,0) -- (3.3,0) node[right] {Real};
       \draw[->] (0,-1.4) -- (0,1.5) node[above] {Imaginary};
    \foreach \x/\xtext in {
0.6/{\scriptstyle\frac{1}{\eta}},
1.13333333/{\scriptstyle{\eta}_1},1.6666666/{\scriptstyle\eta},
3./{{\scriptstyle\eta}+{\scriptstyle\sqrt{{\eta}^2-1}}}}
 \draw[xshift=\x cm] (0pt,1pt) -- (0pt,-1pt) node[below,fill=white]
 {$\xtext$}; 
      \foreach \y/\ytext in {-1.33333/-{\scriptstyle\sqrt{{\eta}^2-1}},
-0.53333/-{\scriptstyle\sqrt{{\eta}_1^2-1}},
      0.533333/{\scriptstyle\sqrt{{\eta}_1^2-1}},
      1.333333/{\scriptstyle\sqrt{{\eta}^2-1}}}
        \draw[yshift=\y cm] (1pt,0pt) -- (-1pt,0pt) node[left,fill=white]
          {$\ytext$};
     \end{scope}
   \draw[arrows=->,style=important line, red] (1.66666666,0) circle (1.3333333333);
   \draw[arrows=->,style=important line, blue] (1.13333333,0) circle (0.5333333333);
 \end{tikzpicture}
\caption{$\mathbb{D}_{\rm inv}(\eta)$~red\quad
$\mathbb{D}_{\rm inv}({\eta}_1)$~blue\quad${\eta}_1
={\scriptstyle\frac{1}{2}\left(\eta+\frac{1}{\eta}\right)}$}
\label{Figure:BasicInvertibleDisks}
\end{figure}

\subsection{Matrix Sign Function Iterations}
\label{Sect:MSF}

In this subsection we review the Matrix Sign Function iteration
scheme, used in matrix computations, see e.g.
\cite[Chapter 5]{Higham2008}, \cite[Chapter 22]{LanRod1995},
\cite{KennLaub1995} and also \cite{CohenLew1997a}, \cite{CohenLew2007},
\cite{Lewk1999a}, \cite{LewRodYar2005}.
\vskip 0.2cm

\noindent
Recall that for a matrix $A\in\C^{n\times n}$ whose spectrum avoids
the imaginary axis, one can define $S:={\rm Sign}(A)\in\C^{n\times n}$
as a matrix satisfying,
\begin{equation}\label{eq:DefSign}
S^2=I_n\quad\quad\quad AS=SA\quad\quad\quad
{\rm spec}(SA)\subset\C_R~.
\end{equation}
The readers interested in learning on the advantages of obtaining the
Sign of a matrix, can look at the above references.
\vskip 0.2cm

\noindent
It is known that for an arbitrary matrix $A\in\C^{n\times n}$ whose
spectrum avoids the imaginary axis, the following iterative
scalar function,
\begin{equation}\label{eq:h_j}
h_o(x)=x\quad{\rm and}\quad
h_j(x)={\scriptstyle\frac{1}{2}}\left(h_{j-1}(x)+{\scriptstyle\frac{1}{h_{j-1}(x)}}
\right)\quad\quad j=1,~2,~\ldots~, 
\end{equation}
is so that 
\[
\lim\limits_{j~\rightarrow~\infty}{h_j(A)}={\rm Sign}(A).
\]
It is also known, see e.g. \cite[Eq. (5.18),~Theorem 5.6]{Higham2008},
that this convergence is quadratic,
for $j\geq 1$, 
\begin{equation}\label{eq:Convergence}
h_j(A)=\left(I-G^{2^j}\right)\left(I+G^{2^j}\right)S
\quad
\quad
\quad
G:=(A-S)(A+S)^{-1},
\end{equation}
We next show how $\mathbb{D}_{\rm inv}(\eta)$ is a natural means of
describing the convergence of the Newton iterations. To this end
we have some preliminaries.\quad
The following is straightforward, yet useful.

\begin{La}\label{La:CayleySignInversion}
Let $\mathcal{C}(A)$ be the Cayley transform of $A\in\C^{n\times n}$ so
that $-1, 0, \pm{i}\not\in{\rm spec}(A)$, then
\begin{itemize}
\item[(i)~~~]{}$-\mathcal{C}\left(A\right)=\mathcal{C}\left(A^{-1}\right)$.
\vskip 0.2cm

\item[(ii)~~]{}$
-\left(\mathcal{C}(A)\right)^2
=\mathcal{C}(\underbrace{{\scriptstyle\frac{1}{2}}(A+A^{-1})
}_{h_1(A)~{\rm see~Eq.}~\eqref{eq:h_j}})$.
\end{itemize}
\end{La}
\vskip 0.2cm

\noindent
Namely, under the Cayley transform $\mathcal{C}$ (see
Definition \ref{Dn:MatrixCayleyTransform}), the operation
of~ {\em inversion} ~takes the form of~ {\em minus}.
\vskip 0.2cm

\noindent
From the above discussion we arrive at the main observation of this
subsection: An illustrative description of the convergence of the Newton
iterations from Eq. \eqref{eq:h_j}, applied to invertible disks.

\begin{Ob}\label{Ob:DisksIterations}
For ${\scriptstyle\eta}\in(1,~\infty)$ let $\mathbb{D}_{\rm
origin}({\scriptstyle\eta})$ and $\mathbb{D}_{\rm inv}(
{\scriptstyle\eta})$ be as in Eqs. \eqref{eq:Dorigin}
and \eqref{eq:InvertibleDisk} respectively. Let also
$h_j$ be as in Eq. \eqref{eq:h_j}.

\noindent
(I) If $A\in\C^{n\times n}$ is so that its spectrum avoids the
imaginary axis, then there exists a finite $\eta>1$ so that
\[
{\rm spec}(A)\subset\{-\mathbb{D}_{\rm inv}(\eta)
\bigcup\mathbb{D}_{\rm inv}(\eta)\}.
\]
Furthermore, for all $j=0,~1,~2,~\ldots$
\[
{\rm spec}\left(h_j(A)\right)\subset\{-\mathbb{D}_{\rm inv}\left(
h_j(\eta)\right)\bigcup\mathbb{D}_{\rm inv}\left(h_j(\eta)\right)\}.
\]
(II)~For $\eta\in(0,~\infty]$ one has that
$\mathcal{C}\left(\mathbb{D}_{\rm inv}(
{\scriptstyle\frac{1}{2}}(\eta+{\scriptstyle\frac{1}{\eta}}))\right)
=-\left(\mathbb{D}_{\rm origin}(\eta)\right)^2$,
and more generally,
\begin{equation}\label{eq:Main}
\mathcal{C}\left(\mathbb{D}_{\rm inv}\left(h_j(\eta)\right)\right)
=(-1)^j\cdot\mathbb{D}_{\rm origin}\left(\left(\sqrt{
{\scriptstyle\frac{\eta-1}{\eta+1}}}\right)^{2^j}\right)
\quad\quad j=0,~1,~2,~\ldots
\end{equation}
\end{Ob}
\vskip 0.2cm

\noindent
Indeed Eq. \eqref{eq:Convergence} and Eq.  \eqref{eq:Main} are quite
similar. An advantage of Observation \ref{Ob:DisksIterations} is that
in each step it provides us with a precise measure $h_j(\eta)$ for
the distance of $h_j(A)$ from $S$.
\vskip 0.2cm

\noindent
This iterative map of the disks of the form
$\mathbb{D}_{\rm inv}(\eta)$, is illustrated in Figure \ref{Figure:MSF}
in stages: From the union of the blue disks, down to the siver disks.
\begin{center}
$
\begin{smallmatrix}
{\eta}=&2+\sqrt{3}&&h_1(\eta)=&2&~&h_2(\eta)=&\frac{5}{4}&~&h_3(\eta)=&
\frac{41}{40}\\~\\
{\rm Radius}_0=&\sqrt{6+4\sqrt{3}}&~&{\rm Radius}_1=&\sqrt{3}&~&
{\rm Radius}_2=&\frac{3}{4}&~&~{\rm Radius}_3=&\frac{9}{40}\\~\\~~~
{\rm blue}&~&&{\rm red}&~&&~~{\rm green}&~&~&~{\rm gray}
\end{smallmatrix}
$
\end{center}
\vskip 0.2cm

\noindent
In Proposition \ref{Pn:Lh(eta)Iterations} below, Observation
\ref{Ob:DisksIterations} is adapted to the framework of Hyper-Lyapunov
inclusions.

\section[Sub-Unit disk]{Matrices whose spectrum is within
a sub-unit disk}
\setcounter{equation}{0}
\label{sec:sub-unitDisk}

\noindent
We now take the scalar parameter $~{\scriptstyle\eta}$ (where
$~{\scriptstyle\eta}\in(1,~\infty]$) from the previous section, and
introduce it to the matricial sets $\mathbf{S}_H$ in Eq.
\eqref{eq:OriginalStein}.
\vskip 0.2cm

\noindent
For arbitrary $H\in\mathbf{H}_n$ and ${\scriptstyle\eta}\in(1,~\infty]$ the
set
\begin{equation}\label{eq:SteinRhoCommonH}
\begin{matrix}
\mathbf{S}_H({\scriptstyle\eta})&=&\left\{ A\in\C^{n\times n}~:~
\left({\scriptstyle\frac{\eta-1}{\eta+1}}H-A^*HA\right)\succ 0~\right\},
\\~\\
\overline{\mathbf S}_H({\scriptstyle\eta})&=&\left\{ A\in\C^{n\times n}~:~
\left({\scriptstyle\frac{\eta-1}{\eta+1}}H-A^*HA\right)\succcurlyeq 0~\right\},
\end{matrix}
\end{equation}
is well defined, e.g. if $U$ is a unitary matrix commuting with
$H$ (i.e. $UH=HU$, $U^*U=I_n$), for all 
\mbox{$\alpha\in(0,~{\scriptstyle\sqrt{\frac{\eta-1}{\eta+1}}})$}
the matrix $A=\alpha U$
belongs to ${\mathbf S}_H({\scriptstyle\eta})$. As before,
$\overline{\mathbf S}_H({\scriptstyle\eta})$ is the closure of the open
set $\mathbf{S}_H({\scriptstyle\eta})$.
\vskip 0.2cm

\noindent
For simplicity, in the rest of this section, we confine the discussion to
$H\succ 0$, i.e. positive definite, and to emphasize that 
we shall denote it by $P$ ($P\succ 0$). As before,
${\scriptstyle\eta}\in(1,~\infty]$  and we shall consider
\begin{equation}\label{eq:SpEta}
\mathbf{S}_P({\scriptstyle\eta})=\left\{ A\in\C^{n\times n}~:~
\left({\scriptstyle\frac{\eta-1}{\eta+1}}P-A^*PA\right)\succ 0~\right\},
\end{equation}
and
$\overline{\mathbf S}_P({\scriptstyle\eta})$.
\vskip 0.2cm

\noindent
Eq. \eqref{eq:SpEta} offers a more detailed examination of
the set $\mathbf{S}_P$, in the following sense,
\[
\infty>\eta>{\eta}_1>1\quad\Longrightarrow\quad
\mathbf{S}_P({\scriptstyle\eta}_1)\subset\mathbf{S}_P({\scriptstyle\eta})
\subset\mathbf{S}_P~,
\]
where each inclusion is strict, and
\[
\lim\limits_{\eta~\longrightarrow~\infty}
\mathbf{S}_P({\scriptstyle\eta})=\mathbf{S}_P~.
\]
Here are three basic properties of this set, which are easy to
verify.

\begin{Pn}\label{Pn:ClosureOfS_H}
For parameters $~P\succ 0$ and \mbox{${\scriptstyle\eta}\in(1,~\infty]$},
the matricial set $\mathbf{S}_P({\scriptstyle\eta})$ in Eq. \eqref{eq:SpEta}
satisfies the following,

\begin{itemize}
\item[(i)~~~]{}For all $c\in\overline{\mathbb D}_{\rm origin}(\eta)$,
see Eq. \eqref{eq:Dorigin}, the set $\mathbf{S}_P({\scriptstyle\eta})$
contains all elements of the form $~cI_n$, and
\[
A\in\mathbf{S}_P({\scriptstyle\eta})\quad
\Longrightarrow\quad
cA\in\mathbf{S}_P({\scriptstyle\eta}).
\]
\item[(ii)~~]{}This set is matrix-product-contractive, i.e. if
$A_a\in\mathbf{S}_P({\scriptstyle{\eta}_a})$ and
$A_b\in\mathbf{S}_P({\scriptstyle{\eta}_b})$
for some
\mbox{${\scriptstyle{\eta}_a}$}, 
\mbox{${\scriptstyle{\eta}_b}\in(1,~\infty]$} then,
the product $A_aA_b$ is s.t.
\[
A_aA_b\in\mathbf{S}_P({\scriptstyle{\eta}_c})
~~{\rm with}~~{\scriptstyle\eta_c}=
{\scriptstyle\frac{1+{\eta}_a{\eta}_b}{{\eta}_a+{\eta}_b}}\quad
{\rm see~Eq.~\eqref{eq:Eta1Eta2}}.
\]
\item[(iii)~]{}Whenever $A_a, A_b\in\mathbf{S}_P({\scriptstyle\eta})$
it implies that
\[
(aA_a+bA_b)\in\mathbf{S}_P\quad\quad\begin{smallmatrix}
\forall~a, b\in\C\\~\\(1+{\scriptstyle\frac{2}{\sqrt{\eta-1}}})
>|a|+|b|.\end{smallmatrix}
\]
\end{itemize}
\end{Pn}

\noindent
By multiplying Eq. \eqref{eq:SpEta} by $P^{-\frac{1}{2}}$ from both sides,
one can equivalently write
\begin{equation}\label{eq:AlternativeSteinRhoCommonH}
\mathbf{S}_P({\scriptstyle\eta})=\{ A\in\C^{n\times n}~:~
{\scriptstyle\sqrt{\frac{\eta-1}{\eta+1}}}
>\| P^{\frac{1}{2}}AP^{-\frac{1}{2}}\|_2~\},
\end{equation}
where $\|~\|_2$ denotes the spectral (a.k.a. Euclidean) norm,
see e.g.  \cite[item 5.6.6]{HornJohnson1}.  
\vskip 0.2cm

\noindent
A a side remark we point out that in \cite[Section 2]{AlpayLew2020}
a rational
functions version of $\overline{\mathbf S}_I({\scriptstyle\eta})$
is studied, i.e.  families of matrix-valued functions $F(z)$ where
\[
{\scriptstyle\sqrt{\frac{\eta-1}{\eta+1}}}
>\| F(z)\|_2\quad\quad\forall z\in\C_R~.
\]
Another set of rational functions will be mentioned in
Subsection \ref{Subsect:Hyper-Lyapunov}
\vskip 0.2cm

\noindent
As a motivation, for studying $\mathbf{S}_P({\scriptstyle\eta})$
one has the following model of a stability robustness problem.
\vskip 0.2cm

\noindent
{\bf Difference inclusion stability problem}\\
Let $x(\cdot)$ be an $n$-dimensional real vector-valued sequence
satisfying,
\begin{equation}\label{eq:DifferenceIncl}
x(k+1)=A\left(k,x(k)\right)x(k)\quad\quad\quad k=0,~1,~2,~\ldots
\end{equation}
where the actual matrix-valued sequence $\{A\left(0,x(0)\right), ~
A\left(1,x(1)\right), ~A\left(2, x(2)\right),~\ldots~\}$ can be~ 
{\em arbitrary}. 
\vskip 0.2cm

\noindent
If in the above difference inclusion,
\[
A(\cdot , \cdot)\in\mathbf{S}_{P}\quad\quad{\rm for~some}~
P\succ 0,
\]
then there exist $\beta$, $\gamma$ where $\beta\geq 1>\gamma\geq 0$
so that,
\[
\beta{\gamma}^k\|x(0)\|_2\geq\|x(k)\|_2\quad\quad\quad
\forall~x(0)\quad k=0,~1,~2,~\ldots
\]
For perspective see e.g. \cite[Eq. (3.3)]{Lewk2020c}, \cite{MolcPyat1989}.
\vskip 0.2cm

\noindent
Note that there is no estimate for $\gamma$, which in principle can be
arbitrarily close to 1.
\vskip 0.2cm

\noindent
On the expense of additional computational burden, more information
is obtained when set $\mathbf{S}_P({\scriptstyle\eta})$ is used.
To simplify the formulation we take $P=I_n~$.

\begin{Pn}\label{Pn:DifferenceInclusion}
If in the above difference inclusion,
\[
A(\cdot , \cdot)\in\mathbf{S}_{I_n}({\scriptstyle\eta}),
\]
for some $\eta\in(1,~\infty)$, then,
\[
\|x(0)\|_2\cdot{\scriptstyle\left(\frac
{\eta-1}{\eta+1}\right)}^{\frac{k}{2}}
\geq\|x(k)\|_2\quad\quad\quad\forall~x(0)\quad
k=0,~1,~2,~\ldots
\]
\end{Pn}
\vskip 0.2cm

\noindent
We end this section by mentioning the following property of
the set $\mathbf{S}_P({\scriptstyle\eta})$.

\begin{Pn}\label{Stein(Eta)}
For $P=I$ and arbitrary \mbox{${\scriptstyle\eta}\in(1, \infty]$}
the matricial set $\mathbf{S}_I({\scriptstyle\eta})$ in Eqs.
\eqref{eq:SteinRhoCommonH} and \eqref{eq:AlternativeSteinRhoCommonH},
is matrix-convex.
\end{Pn}
\vskip 0.2cm

\noindent
This can be deduced either: From part (I) of \cite[Observation 2.2]{Lewk2020a},
or: From Lemma \ref{La:MatrixConvex} below.

\section{Hyper-Lyapunov Inclusions}
\setcounter{equation}{0}

\subsection{Hyper-Lyapunov Inclusion - Properties}
\label{Subsect:Hyper-Lyapunov}

\noindent
We start with a basic fact, whose verification is straightforward.

\begin{Pn}\label{Pn:CayleyLyapStein}
For $H\in\mathbf{H}_n$ and ${\scriptstyle\eta}\in(1,~\infty]$, both
arbitrary, let $\mathbf{S}_H({\scriptstyle\eta})$ be as in Eq.
\eqref{eq:SteinRhoCommonH}. Then, one has that
\begin{equation}\label{eq:CayleySteinEta}
\mathcal{C}\left(\mathbf{S}_H({\scriptstyle\eta})\right)
=\mathbf{L}_H({\scriptstyle\eta}),
\end{equation}
where,
\begin{equation}\label{eq:DefLh(eta)}
\mathbf{L}_H({\scriptstyle\eta})=\left\{ A\in\C^{n\times n}~:~
(HA+A^*H)\succ{\scriptstyle\frac{1}{\eta}}(A^*HA+H)~\right\}.
\end{equation}
\end{Pn}

\begin{Rk}\label{Rk:SteinLyapRhoCommonH}
{\rm
Note that in the {\em Hyper-Stein} inclusion \eqref{eq:SteinRhoCommonH},
\eqref{eq:CayleySteinEta} and in the {\em Hyper-Lyapunov} inclusion
\eqref{eq:DefLh(eta)}, the parameters $H$ and ${\scriptstyle\eta}$,
are indeed the same.
}
\qed
\end{Rk}

\noindent
Here is the first result concerning this set.

\begin{Pn}\label{Pn:FirstProperties}
For arbitrary $H\in\mathbf{H}_n$ let ${\rm Sign}(H)$ be as defined in
Eq. \eqref{eq:DefSign} and let ${\scriptstyle\eta}\in(1,~\infty]$,

\begin{itemize}

\item[(i)~~~]{}The set $\mathbf{L}_H({\scriptstyle\eta})$ is
closed under inversion.
\vskip 0.2cm 

\item[(ii)~~]{}The matrix $A=c\cdot{\rm Sign}(H)$ belongs to 
\mbox{$\mathbf{L}_H({\scriptstyle\eta})$} for all
\mbox{$c\in\overline{\mathbb D}_{\rm inv}({\scriptstyle\eta})$, see Eq.
\eqref{eq:InvertibleDisk}}.
\end{itemize}
\end{Pn}

\noindent
{\bf Proof :}~ (i)~~We shall show inveribility in two ways:\\
(a) Multiply Eq. \eqref{eq:DefLh(eta)} by $\left(A^{-1}\right)^*$
and $A^{-1}$ from the left and from the right, respectively. Then,
\[
\begin{matrix}
~&\left(A^{-1}\right)^*\left((HA+A^*H)\succ\frac{1}{\eta}(A^*HA+H)
\right)A^{-1}
\\
{\rm results~in}&~\\
~&(HA^{-1}+(A^{-1})^*H)\succ\frac{1}{\eta}((A^{-1})^*HA^{-1}+H).
\end{matrix}
\]
(b)~~Invertibility can be indirectly deduced by using the relation
$\mathcal{C}\left(A^{-1}\right)=-\mathcal{C}\left(A\right)$ 
which is immediate from Definition \ref{Dn:MatrixCayleyTransform},
along with Eqs. \eqref{eq:SteinRhoCommonH} and \eqref{eq:DefLh(eta)}.
\vskip 0.2cm

\noindent
(ii) Substitute $A=c{\rm Sign}(H)$ into each side of Eq.
\eqref{eq:DefLh(eta)} then on the left-hand side:
\[
\begin{matrix}
\eta(HA+A^*H)=2\eta\cdot{\rm Re}(c)H{\rm Sign}(H)\succcurlyeq(1+|c|^2)
H{\rm Sign}(H)&&{\rm by~item}~(i)~{\rm of~Lemma~\ref{La:Dinv}},
\end{matrix}
\]
and on the right-hand side $~ A^*HA+H=(1+|c|^2)H$.\\
As always, $H{\rm Sign}(H)\succ H$, so this claim is established.
\qed
\vskip 0.2cm 

\noindent
At this point we find it convenient to confine the discussion in the
sequel to having $H=P$ positive definite ($H$ suggests ``Hermitian"
while $P$ stands for``positive"), so Eq. \eqref{eq:DefLh(eta)}
takes the form
\begin{equation}\label{eq:DefLp(eta)}
\mathbf{L}_P({\scriptstyle\eta})=\left\{ A\in\C^{n\times n}~:~
(PA+A^*P)\succ{\scriptstyle\frac{1}{\eta}}(A^*PA+P)~\right\},
\end{equation}
where as before, ${\scriptstyle\eta}\in(1,~\infty]$.
\quad
Expectedly, in Eq. \eqref{eq:DefLp(eta)} we have
\[
\begin{matrix}
\lim\limits_{\eta~\longrightarrow~\infty}\mathbf{L}_P({\scriptstyle\eta})
&=&
\lim\limits_{\eta~\longrightarrow~\infty}\left\{ A\in\C^{n\times n}~:~
(PA+A^*P)\succ{\scriptstyle\frac{1}{\eta}}(A^*PA+P)~\right\}
\\~\\~&=&
\left\{ A\in\C^{n\times n}~:~\left(PA+A^*P\right)\succ 0~\right\}
=\mathbf{L}_P,
\end{matrix}
\]
where the lower part is as in Eq. \eqref{eq:OriginalLyapunov}.
Furthermore, for a prescribed $~P\succ 0$ one has that,
\begin{equation}\label{eq:Inclusions}
\infty>{\scriptstyle\eta}>{\scriptstyle{\eta}_1}>1\quad\Longrightarrow
\quad\mathbf{L}_P({\scriptstyle{\eta}_1})\subset
\mathbf{L}_P({\scriptstyle\eta})\subset\mathbf{L}_P~,
\end{equation}
and each inclusion is strict.
\vskip 0.2cm 

\noindent
This family can also be written as,
\[
\mathbf{L}_P({\scriptstyle\eta})=\left\{ A\in\C^{n\times n}~:~
(A+I_n)^*P(A+I_n)\succ(1+{\scriptstyle\frac{1}{\eta}})
(A^*PA+P)~\right\},
\]
which illustrates the fact that $\mathbf{L}_P({\scriptstyle\eta})$
is (linear in $P$ and) quadratic in $A$ (see also Subsection
\ref{SubSec:UnifyingFramework}), so this set may contain $A$ only
when $\| A\|$ is neither ``too small" nor ``too large".
\vskip 0.2cm 

\noindent
We next present yet another description of the set
$\mathbf{L}_P({\scriptstyle\eta})$. To this end, denote
by $\rho(M)$ the spectral radius of a square matrix $M$.

\begin{Cy}\label{Cy:QuantitativeHyperPositive1}
For $P\succ 0$ and ${\scriptstyle\eta}\in(1,~\infty]$ let
$\mathbf{L}_P({\scriptstyle\eta})$ be as before, Then,
the following is true.
\begin{itemize}
\item[(i)~~~]{}
$~\mathbf{L}_P({\scriptstyle\eta})=\left\{~A\in\mathbf{L}_P~:~
{\scriptstyle\eta}>\rho\left({\scriptstyle
\left(A^*PA+P\right)\left(A^*P+PA\right)^{-1}}\right)\right\}.$
\vskip 0.2cm

\item[(ii)~~]{}$~A\in\mathbf{L}_P(\eta)$ implies that
${\rm spec}(A)\subset\mathbb{D}_{\rm inv}(\eta)$.
\vskip 0.2cm

\item[(iii)~]{} $~A\in\mathbf{L}_P(\eta)~$ implies that for all
\mbox{$\gamma\in[{\scriptstyle 1,~\frac{\eta}{\sqrt{{\eta}^2
-1}}})$} the matrix $A^{\gamma}$ belongs to $\mathbf{L}_P~$.
\end{itemize}
\end{Cy}

\noindent
{\bf Proof :}\quad
(i)~Writing Eq. \eqref{eq:DefLp(eta)} explicitly, one has the following
chain of relations
\[
\begin{matrix}
\mathbf{L}_P({\scriptstyle\eta})&=&\left\{ A\in\C^{n\times n}~:\right.&
\left.(PA+A^*P)\succ{\scriptstyle\frac{1}{\eta}}(A^*PA+P)~\right\},
\\~\\~
&=&\left\{~A\in\mathbf{L}_P~:\right.&\left.
(PA+A^*P)\succ{\scriptstyle\frac{1}{\eta}}(A^*PA+P)~\right\},
\\~\\~
&=&\left\{~A\in\mathbf{L}_P~:\right.&\left.
I_n\succ{\scriptstyle\frac{1}{\eta}}(PA+A^*P)^{-\frac{1}{2}}
(A^*PA+P)(PA+A^*P)^{-\frac{1}{2}}~\right\},
\\~\\~
&=&\left\{~A\in\mathbf{L}_P~:\right.&\left.
\eta>\|(PA+A^*P)^{-\frac{1}{2}}
(A^*PA+P)(PA+A^*P)^{-\frac{1}{2}}\|_2\right\},
\\~\\~
&=&\left\{~A\in\mathbf{L}_P~:\right.&\left.
\eta>\rho\left(~(PA+A^*P)^{-\frac{1}{2}}
(A^*PA+P)(PA+A^*P)^{-\frac{1}{2}}\right)~\right\},
\\~\\~
&=&\left\{~A\in\mathbf{L}_P~:\right.&\left.
\eta>\rho\left(~(A^*PA+P)(PA+A^*P)^{-1}\right)~
\right\},
\end{matrix}
\]
so this part of the claim is established.
\vskip 0.2cm

\noindent
(ii) By Eq. \eqref{eq:AlternativeSteinRhoCommonH}
\[
B\in\mathbf{S}_P(\eta)\quad\Longrightarrow\quad
{\rm spec}(B)\subset\mathbb{D}_{\rm origin}(
{\scriptstyle\sqrt{\frac{\eta-1}{\eta+1}}}),
\]
and thus applying the Cayley transform
\[
\underbrace{\mathcal{C}(B)}_{A}\in\underbrace{\mathcal{C}(
\mathbf{S}_P(\eta))}_{\mathbf{L}_P(\eta),~{\rm by~Eq.}~
\eqref{eq:CayleySteinEta}}\quad\Longrightarrow\quad
{\rm spec}(\underbrace{\mathcal{C}(B)}_{A})\subset
\underbrace{\mathcal{C}(\mathbb{D}_{\rm origin}(
{\scriptstyle\sqrt{\frac{\eta-1}{\eta+1}}}))}
_{\mathbb{D}_{\rm inv}(\eta),~{\rm by~Lemma}~\ref{La:Dinv}~{\rm (iv)}},
\]
yields item (ii) of the claim.
\vskip 0.2cm

\noindent
(iii) From the previous item we know that
${\rm spect}(A)\subset\mathbb{D}_{\rm inv}(\eta)$. By definition,
$\mathbb{D}_{\rm inv}(\eta)$ is a disk of radius
${\scriptstyle\sqrt{{\eta}^2-1}}$, centered at
${\scriptstyle\eta+i0}$, see Eq. \eqref{eq:InvertibleDisk}. Hence,
$\mathbb{D}_{\rm inv}(\eta)$ intersects the unit disk at the points
${\scriptstyle\frac{1}{\eta}\pm{i}\frac{\sqrt{{\eta}^2-1}}{\eta}}$,
see Figure \ref{Figure:Eta=sqrt{2}}. Thus, the rays emanating from
the origin, passing through these points, are tangential to 
$\mathbb{D}_{\rm inv}(\eta)$. One can conclude that
$\left(\mathbb{D}_{\rm inv}(\eta)\right)^{\gamma}\subset\C_R$
for all prescribed $\gamma$, so the claim is established
\qed
\vskip 0.2cm

\noindent
Note that having $A\in\mathbf{L}_P$, for some $P\succ 0$, does not
provide us
with any information on the distribution of the eigenvalues of $A$,
within $\C_R~$. Items (ii) and (iii) of the above Corollary
\ref{Cy:QuantitativeHyperPositive1}
indicates that ${\rm spec}(A)\subset{\mathbb D}_{\rm inv}(\eta)$,
see e.g. Figure \ref{Figure:Eta=sqrt{2}}.
Needless to say, this information on the spectrum comes on the
expense of an additional computational burden. 
\vskip 0.2cm

\noindent
We next illustrate the above results.

\begin{Ex}
{\rm
For a parameter ${\scriptstyle\theta}\in[0,~2\pi)$ consider a
pair of structured real $2\times 2$ matrices, along with the
corresponding Cayley transform and the respective spectrum,
\[
\begin{matrix}
A_1={\scriptstyle\eta}I_2
-{\scriptstyle\sqrt{{\eta}^2-1}}
\left(\begin{smallmatrix}\cos(\theta)
&&~\sin(\theta)
\\~\\
\sin(\theta)
&&-\cos(\theta)
\end{smallmatrix}\right)
&&&
\mathcal{C}(A_1)
={\scriptstyle\frac{\sqrt{\eta-1}}{\sqrt{\eta+1}}}
\left(\begin{smallmatrix}\cos(\theta)&&~\sin(\theta)\\~\\
\sin(\theta)&&-\cos(\theta)\end{smallmatrix}\right)
\\~\\
{\scriptstyle{\lambda}_{1,2}}(A_1)={\scriptstyle\eta\pm\sqrt{{\eta}^2-1}}
&&&
{\scriptstyle{\lambda}_{1,2}}\left(\mathcal{C}(A_1)\right)
=\pm{\scriptstyle\frac{\sqrt{\eta-1}}{\sqrt{\eta+1}}}
\end{matrix}
\]
and
\[
\begin{matrix}
A_2=
{\scriptstyle\frac{1}
{\eta+\cos(\theta)\sqrt{{\eta}^2-1}}}
\left(\begin{smallmatrix}1&&
\sin(\theta)\sqrt{{\eta}^2-1}
\\~\\
-\sin(\theta)\sqrt{{\eta}^2-1}&&1\end{smallmatrix}\right)
&&&
\mathcal{C}(A_2)
={\scriptstyle\frac{\sqrt{\eta-1}}{\sqrt{\eta+1}}}
\left(\begin{smallmatrix}\cos(\theta)&&-\sin(\theta)\\~\\
\sin(\theta)&&~\cos(\theta)\end{smallmatrix}\right)
\\~\\
{\scriptstyle{\lambda}_{1,2}}(A_2)=
{\scriptstyle\frac{1\pm{i}\sin(\theta)\sqrt{{\eta}^2-1}}
{\eta+\cos(\theta)\sqrt{{\eta}^2-1}}}
&&&
{\scriptstyle{\lambda}_{1,2}}\left(\mathcal{C}(A_2)\right)
={\scriptstyle\frac{\sqrt{\eta-1}}{\sqrt{\eta+1}}}\left(
{\scriptstyle\cos(\theta)\pm{i}\sin(\theta)}\right).
\end{matrix}
\]
Indeed, as $\mathcal{C}(A_1)$ and $\mathcal{C}(A_2)$ are
scaled orthogonal matrices, their spectrum is on the boundary
of $\overline{\mathbb D}_{\rm origin}\left({\scriptstyle\frac
{\sqrt{\eta-1}}{\sqrt{\eta+1}}}\right)$; and the spectrum of
$A_1$ and $A_2$ is on the boundary of
$\overline{\mathbb D}_{\rm inv}(\eta)$.\\
Furthermore, here item (i) of Corollary
\ref{Cy:QuantitativeHyperPositive1} takes a special form,

\begin{center}
$(A_1^*A_1+I_2)(A_1+A_1^*)^{-1}=\eta{I}_2
=(A_2^*A_2+I_2)(A_2+A_2^*)^{-1}.$
\end{center}
}
\qed
\end{Ex}
\vskip 0.2cm

\noindent
We now show the main properties of $\mathbf{L}_H(\eta)$, cited
in the Introduction.

\noindent
{\bf Proof of Theoren \ref{Tm:ClosureOfL_H(r)} :}\quad  (i)
Invertibility was shown in item (i) of
Proposition \ref{Pn:FirstProperties}.
\vskip 0.2cm

\noindent
To show convexity, for a pair of $n\times n$ matrices
$A_o$, $A_1$ consider the following identity, 
\[
\begin{matrix}
{\scriptstyle\theta}
\underbrace{
\left(
-{\scriptstyle\frac{1}{\eta}}({A_1}^*PA_1+P)+PA_1+{A_1}^*P
\right)}_{Q_1}
\\~\\
+
\left({\scriptstyle 1}-{\scriptstyle\theta}\right)
\underbrace{\left(-{\scriptstyle\frac{1}{\eta}}({A_o}^*+PA_o+P)
+PA_o+{A_o}^*P\right)}_{Q_o}
\\~\\
+
{\scriptstyle\frac{\theta(1-\theta)}{\eta}}\left(A_o-A_1\right)^*P\left(
A_o-A_1\right)
\\~\\=
\left(-{\scriptstyle\frac{1}{\eta}}({A_{\theta}}^*PA_{\theta}+P)
+PA_{\theta}+{A_{\theta}}^*P\right)\quad\quad\quad
A_{\theta}:=
{\scriptstyle\theta}{A}_1+{\scriptstyle(1-\theta)}A_o~.
\end{matrix}
\]
The condition $~A_o, A_1\in\mathbf{L}_P({\scriptstyle\eta})$ is equivalent
to having $Q_o, Q_1\succ 0$. Since by assumption  $P\succ 0$, it implies
that \mbox{$\left(A_o-A_1\right)^*P\left(
A_o-A_1\right)\succcurlyeq 0$} and thus one has
that for all \mbox{${\scriptstyle\theta}\in[0,~1]$}, also 
\mbox{$A_{\theta}:=\left({\scriptstyle\theta}{A}_1+
(1-{\scriptstyle\theta})A_o\right)$}
belongs to $~\mathbf{L}_P({\scriptstyle\eta})$,
so this part of the claim is established.
\vskip 0.2cm

\noindent
(ii)~~ This is shown in two ways: (a) Combine item (i) along with
Eq. \eqref{eq:Inclusions}; (b) See e.g.
\cite[Lemma 3.5]{CohenLew1997a}.
\vskip 0.2cm

\noindent
(iii)~~Matrix convexity for $P=I$: In Lemma \ref{La:MatrixConvex}
below, substitute $\alpha=\gamma=-{\scriptstyle\frac{1}{\eta}}$
and $\beta=1$. 
\vskip 0.2cm

\noindent
(iv)~~Multiply Eq. \eqref{eq:DefLp(eta)} by $T^*$ and 
$T$ from the left and from the right respectively,
\[
\begin{matrix}
T^*\left(~(PA+A^*P)\succ{\scriptstyle\frac{1}{\eta}}(A^*PA+P)~\right)T
\\~\\
\underbrace{T^*PT}_{\hat{P}}\underbrace{T^{-1}AT}_{\hat{A}}+
\underbrace{T^*A^*(T^{-1})^*}_{\hat{A}^*}\underbrace{T^*PT}_{\hat{P}}
\succ{\scriptstyle\frac{1}{\eta}}
\underbrace{T^*A^*(T^{-1})^*}_{\hat{A}^*}\underbrace{T^*PT}_{\hat{P}}
\underbrace{T^{-1}AT}_{\hat{A}}+\underbrace{T^*PT}_{\hat{P}}.
\end{matrix}
\]
Thus, the proof is complete.
\qed
\vskip 0.2cm

\noindent
A a side remark we point out that in \cite[Section 3]{AlpayLew2020}
a rational functions version of
$\overline{\mathbf L}_I({\scriptstyle\eta})$ is studies, i.e.
families of $m\times m$-valued functions $F(z)$ where
\[
\eta\geq\sup\limits_{z\in\C_R}\rho\left(
(F^*F+I_m)(F^*+F)^{-1}\right).
\]
It is subsequently showed that these are used in the context of
absolute stability of feedback-loops,
for background see e.g. \cite[Section 7.1]{Khalil2000}
\vskip 0.2cm

\noindent
In the previous section we applied the set $\mathbf{S}_P(\eta)$ to
the difference inclusion problem. To simplify the exposition,
instead of taking the continuous-time analogue of differential
inclusion, we here  address the question of uncertain linear,
autonomous system.
\vskip 0.2cm

\noindent
{\bf Stability of uncertain linear autonomous system}\\
For a given set $\mathbf{A}$ of real $n\times n$ matrices
consider the equation,
\begin{equation}\label{eq:UncertainLTI}
\dot{x}=Ax\quad\quad\quad A\in\mathbf{A},
\end{equation}

\hfill{where $A$ is not prescribed.}
\vskip 0.2cm

\noindent
If \mbox{$\mathbf{A}\subset\mathbf{L}_{-P}$}, for some $P\succ 0$, 
then there exist $\beta\geq 1$ and $\alpha>0$ so that
\[
\beta e^{-t\alpha}\|x(0)\|_2\geq\|x(k)\|_2
\quad\quad\quad\forall~x(0)\quad t\geq 0.
\]
For perspective see e.g. 
\cite[Observation 2.1]{Lewk2020a}, \cite{MolcPyat1989}.
\vskip 0.2cm

\noindent
Note that there is no estimate of $\alpha$, nor on the distribution
of trajectories $x(t)$. Namely, with the same $A$, for two distinct
initial conditions, say $x_a(0)$, $x_b(0)$ which are (normalized)
eigenvectors of $A$, one may have,
\[
x_a(t)=e^{-\alpha{t}\pm\gamma\sin(2\pi\alpha{t})}x_a(0),
\quad\quad
x_b(t)=e^{-\delta\alpha{t}}x_b(0),
\quad 
\begin{smallmatrix}
\gamma, \delta>>1\\~\\
t\geq 0.
\end{smallmatrix}
\]
In particular, this implies that $\frac{|\log(\|x_b(t)\|)|}{|\log(\|x_a(t)\|)|}$
may be ``large".
\vskip 0.2cm

\noindent
We next show that by substituting $\mathbf{L}_P$ with
\mbox{$\mathbf{L}_P({\scriptstyle\eta})$}, the above observation can
be refined to simultaneously {\em cluster} all trajectories.
\vskip 0.2cm

\noindent
Recall, see Figures \ref{Figure:BasicInvertibleDisks} and
\ref{Figure:Eta=sqrt{2}}, that a point within
$-\mathbb{D}_{\rm inv}(\eta)$, can be parametrized as
\[
-\eta+\rho\cos(\phi)+i\rho\sin(\phi)\quad\quad\quad
\begin{smallmatrix}
\rho\in\left[0,~\sqrt{{\eta}^2-1}~\right)\\~\\
\phi\in[0,~2\pi).
\end{smallmatrix}
\]
Consider the problem in Eq. \eqref{eq:UncertainLTI}, where now
\mbox{$\mathbf{A}\subset\mathbf{L}_{-P}(\eta)$}, for some
$P\succ 0$ and some $\eta\in(1,~\infty]$. Then for $k=1, 2$
whenever $x_1(0)$ $x_2$ are eigenvectors of $A_1$,
$A_2$ both in
\mbox{$\mathbf{L}_{-P}({\scriptstyle\eta})$},
there are ${\rho}_k, {\phi}_k, {\theta}_k$ so that,
\[
x_k(t)=x_k(0)e^{t(-\eta+{\rho}_k\cos({\phi}_k))}\sin\left(
{\scriptstyle{\theta}_k}+({\rho}_k\sin({\phi}_k))t\right)
\quad\quad\quad
\begin{smallmatrix}
k=1,~2&&{\rho}_k\in\left[0,~\sqrt{{\eta}^2-1}~\right)\\~\\
{\phi}_k, {\theta}_k\in[0,~2\pi)&&\forall t\geq 0.
\end{smallmatrix}
\]
One can now conclude that any two trajectories $x_1(t)$ and $x_2(t)$,
must be in the vicinity of each other in the sense that
$\frac{|\log(\|x_2(t)\|)|}{|\log(\|x_1(t)\|)|}$ is ``small" for
all $t\geq 0$.
\vskip 0.2cm

\noindent
We end this subsection by further relating the new concept of
Hyper-Lyapunov inclusion in Eq. \eqref{eq:DefLh(eta)} to the
classical notion of Matrix Sign Function iterations, from Section
\ref{Sect:MSF}, see e.g. \cite[Chapter 5]{Higham2008},
\cite[Chapter 22]{LanRod1995}, \cite{KennLaub1995} and also
\cite{CohenLew1997a}, \cite{CohenLew2007}, \cite{Lewk1999a},
\cite{LewRodYar2005}.

\begin{Pn}\label{Pn:Lh(eta)Iterations}
Let $A$ be a matrix whose spectrum is in $\C_R$ (a.k.a. ``positively
stable"). Then, there exist parameters $\eta\in(1,~\infty)$ and
(non-unique) $P\succ 0$ so that
\mbox{$A\in\mathbf{L}_P({\scriptstyle\eta})$}. Moreover, the matrix
${\scriptstyle\frac{1}{2}}(A+A^{-1})$ is well defined and
(with the same $P$ and $\eta$)
belongs to \mbox{$\mathbf{L}_P\left({\scriptstyle\frac{1}{2}}\left(
{\scriptstyle\eta+\frac{1}{\eta}}\right)\right)$}.
\end{Pn}

\noindent
{\bf Proof :}~ For a positively stable $A$, it is well known, see
e.g. \cite[Theorem 2.2.3]{HornJohnson2} that for all $Q\succ 0$
there exists $P\succ 0$ so that $PA+A^*P=Q$. Then the following
spectral radius is well defined
\mbox{$\rho\left((A^*PA+P)(PA+A^*P)^{-1}\right)$}, and by 
item (i) of Corollary \ref{Cy:QuantitativeHyperPositive1}
indeed \mbox{$A\in\mathbf{L}_P({\scriptstyle\eta})$}.
\vskip 0.2cm

\noindent
It is also well known, see e.g. \cite[Chapter 5]{Higham2008},
\cite[Chapter 22]{LanRod1995}, \cite{KennLaub1995} and also
\cite{CohenLew1997a}, \cite{CohenLew2007}, \cite{Lewk1999a},
\cite{LewRodYar2005}
that \mbox{$A_1:={\scriptstyle\frac{1}{2}}(A+A^{-1})$} is well
defined and by item (i) of Theorem \ref{Tm:ClosureOfL_H(r)}
one has that $A_1\in\mathbf{L}_P({\scriptstyle\eta})$.
\vskip 0.2cm

\noindent
To show that $A_1$ in fact belongs to
\mbox{$\mathbf{L}_P\left({\scriptstyle\frac{1}{2}}\left(
{\scriptstyle\eta+\frac{1}{\eta}}\right)\right)$}, we return to the Cayley
transform and recall that from Eq. \eqref{Pn:CayleyLyapStein} that
\mbox{$\mathcal{C}(A)\in\mathbf{S}_P({\scriptstyle\eta})$} and by items
(ii) and (iii) of Proposition \ref{Pn:ClosureOfS_H} it follows that
$-\left(\mathcal{C}(A)\right)^2$ belongs to
$\mathbf{S}_P\left({\scriptstyle\frac{1}{2}}\left(
{\scriptstyle\eta+\frac{1}{\eta}}\right)\right)$. Applying again the Cayley
transform, by item (ii) of Lemma \ref{La:CayleySignInversion}, the
construction is complete.
\qed
\vskip 0.2cm

\noindent
See again Figures \ref{Figure:BasicInvertibleDisks} and
\ref{Figure:MSF}.
\vskip 0.2cm

\noindent
Proposition \ref{Pn:Lh(eta)Iterations} can be formulated for the case
where \mbox{$A_1={\scriptstyle\theta}A+({\scriptstyle 1-\theta})A^{-1}$},
for arbitrary ${\scriptstyle\theta}\in[0,~1]$. The choice
\mbox{${\scriptstyle\theta}={\scriptstyle\frac{1}{2}}$} was
adopted just to simplify the exposition. For details, see items (i)
and (iii) of Lemma \ref{La:ProductHBeta}.
\vskip 0.2cm

\subsection{Quadratic Inclusions - a Unifying Framework}
\label{SubSec:UnifyingFramework}

We now introduce a unifying framework for the four inclusions presented
thus far. To this end let $M\in\mathbf{H}_{n+m}$ be a
matricial parameter (i.e. $(n+m)^2$ real scalar parameters)
\begin{equation}\label{eq:BasicRelation}
\left(\begin{smallmatrix}A\\~\\I_m\end{smallmatrix}\right)^*
\underbrace{\left(\begin{smallmatrix}W&&{R}^*\\~\\R&&Y
\end{smallmatrix}\right)}_{M\in\mathbf{H}_{n+m}}
\left(\begin{smallmatrix}A\\~\\I_m\end{smallmatrix}\right)
=\underbrace{\begin{smallmatrix}A^*WA+RA+(RA)^*+Y
\end{smallmatrix}}_{Q=Q^*~m\times m},
\end{equation}
i.e. $A\in\C^{n\times m}$ where $n$ may be larger, equal or smaller
than $m$.
\vskip 0.2cm

\noindent
A word of caution. Our notation, aimed at being consistent throughout
the work, is not common in the Riccati equation circles.  Nevertheless,
there is no contradiction: When in Eq. \eqref{eq:BasicRelation} $m=n$,
it does conform with conventional Riccati equation, see e.g.
\cite{LanRod1995},
\[
\left(\begin{smallmatrix}I_n&-A^*\end{smallmatrix}\right)
\underbrace{
\left(\begin{smallmatrix}R&&~Y\\~\\W&&-R^*\end{smallmatrix}\right)
}_{\rm Hamiltonian}
\left(\begin{smallmatrix}A\\~\\I_n\end{smallmatrix}\right)=
0_{n\times n}
\]
As a specific example, the classical Linear Quadratic Regulator control
problem is obtained when in Eq. \eqref{eq:BasicRelation} one takes
$Q=0$, $Y\succcurlyeq 0$ so that the pair $R, Y$ is detectable and
$W=-B\hat{W}^{-1}B$, where $B\in\C^{n\times m}$ is so that the pair
$R, B$ is stabilizable and $\hat{W}\succ 0$, see e.g.
\cite[Section 6.2]{DullPaga2000}, \cite[Chapter 16]{LanRod1995}.
\vskip 0.2cm

\noindent
We next focus on the case where in Eq. \eqref{eq:BasicRelation}
$n=m$, and Eq. \eqref{eq:BasicRelation} takes the form,
\begin{equation}\label{eq:QuadraticForm}
\left(\begin{smallmatrix}A\\~\\I_n\end{smallmatrix}\right)^*
M\left(\begin{smallmatrix}A\\~\\I_n\end{smallmatrix}\right)
\succcurlyeq 0.
\end{equation}
More specifically, with parameters $~H\in\mathbf{H}_n$,
$P\succ 0$ and $\eta\in(1,~{\infty}]$,
the four previous cases can be recovered,
\[
\begin{array}{c|l}
M&{\rm Inclusion~name}\\
\hline
~&~\\
\left(\begin{smallmatrix}0&&H\\~\\H&&0\end{smallmatrix}\right)
&~{\rm Lyapunov,~~~Eq.}~\eqref{eq:OriginalLyapunov}\\
~&~\\
\left(\begin{smallmatrix}-H&&0\\~\\~0&&H\end{smallmatrix}\right)
&~{\rm Stein,~~~Eq.}~\eqref{eq:OriginalStein}\\
~&~\\
\left(\begin{smallmatrix}-P&&0\\~\\~0&&
\frac{\eta-1}{\eta+1}P\end{smallmatrix}\right)
&~{\rm Hyper-Stein,~~~Eq.}~\eqref{eq:SteinRhoCommonH}\\
~&~\\
\left(\begin{smallmatrix}
-\frac{1}{\eta}P&&~~P\\~\\~P&&-\frac{1}{\eta}P\end{smallmatrix}\right)
&~{\rm Hyper-Lyapunov,~~~Eq.}~\eqref{eq:DefLh(eta)}.
\end{array}
\]

\begin{Rk}
{\rm 
In each of the four above cases, the matrix $M$ has $n$ positive and
$n$ negative eigenvalues (since by assumption, $H\in\mathbf{H}_n$).
In fact, in Eq. \eqref{eq:QuadraticForm} the $2n\times n$ matrix
$\left(\begin{smallmatrix}A\\~\\I_n\end{smallmatrix}\right)$
corresponds to a non-negative subspace of $M$.
\qed
}
\end{Rk}
\vskip 0.2cm

\noindent
In \cite{Lewk2020b} the idea of the above table is extended to
four variants of state-space realization arrays.
\vskip 0.2cm

\noindent
The following technical result is useful in showing matrix-convexity
of various sets.

\begin{La}\label{La:MatrixConvex}
Let the real parameters $\alpha$, $\beta$ and $\gamma$, be so that the
set of matrices $A\in\C^{n\times n}$ satisfying,
\begin{equation}\label{eq:BasicMatrixConvex}
\left(\begin{smallmatrix}A\\~\\I_n\end{smallmatrix}\right)^*
\left(\begin{smallmatrix}
{\alpha}I_n&&{\beta}I_n\\~\\{\beta}I_n&&{\gamma}I_n
\end{smallmatrix}\right)
\left(\begin{smallmatrix}A\\~\\ I_n\end{smallmatrix}\right)
\succcurlyeq 0,
\end{equation}
is not empty. Then, the family of all matrices $A$ satisfying that,
is matrix-convex.
\end{La}

\noindent
{\bf Proof :}~Indeed for $j=1,~\ldots~,~k$ let $A_j$ be so that, 
\[
\left(\begin{smallmatrix}A_j\\~\\I_n
\end{smallmatrix}\right)^*
\left(\begin{smallmatrix}{\alpha}I_n&&{\beta}I_n\\~\\{\beta}I_n&&
{\gamma}I_n\end{smallmatrix}\right)\left(\begin{smallmatrix}A_j
\\~\\I_n\end{smallmatrix}\right)=\underbrace{\begin{smallmatrix}
\alpha{A_j}^*A_j+\beta(A_j+{A_j}^*)+{\gamma}I_n\end{smallmatrix}
}_{Q_j}\succcurlyeq 0.
\]
Now, by Definition \ref{Dn:MatrixConvex}, for ${\upsilon}_j$ so
that $\sum\limits_{j=1}^k{{\upsilon}_j}^*{\upsilon}_j=I_n$
one has that,
\[
\begin{matrix}
\left(\begin{smallmatrix}
\sum\limits_{j=1}^k{{\upsilon}_j}^*A_j{\upsilon}_j
\\~\\I_n\end{smallmatrix}\right)^*\left(\begin{smallmatrix}{\alpha}I_n
&&{\beta}I_n\\~\\{\beta}I_n&&{\gamma}I_n\end{smallmatrix}\right)
\left(\begin{smallmatrix}
\sum\limits_{j=1}^k{{\upsilon}_j}^*A_j{\upsilon}_j
\\~\\ I_n\end{smallmatrix}\right)
=
\left(\begin{smallmatrix}
\sum\limits_{j=1}^k{{\upsilon}_j}^*A_j{\upsilon}_j
\\~\\
\sum\limits_{j=1}^k{\upsilon}_j^*{\upsilon}_j
\end{smallmatrix}\right)^*
\left(\begin{smallmatrix}
{\alpha}I_n&&{\beta}I_n\\~\\{\beta}I_n&&{\gamma}I_n
\end{smallmatrix}\right)
\left(\begin{smallmatrix}
\sum\limits_{j=1}^k{{\upsilon}_j}^*A_j{\upsilon}_j
\\~\\
\sum\limits_{j=1}^k{{\upsilon}_j}^*{\upsilon}_j
\end{smallmatrix}\right)
\\~\\=
\begin{smallmatrix}
\sum\limits_{j=1}^k{{\upsilon}}_j^*\left(
{\alpha}A_j^*A_j+{\beta}({A_j}^*+A_j)+{\gamma}I_n\right){\upsilon}_j
\end{smallmatrix}
=
\sum\limits_{j=1}^k{{\upsilon}}_j^*
\underbrace{
\left(\begin{smallmatrix}{A_j}^*&&I_n\end{smallmatrix}\right)
\left(\begin{smallmatrix}{\alpha}I_n&&{\beta}I_n\\~\\{\beta}I_n&&
{\gamma}I_n\end{smallmatrix}\right)
\left(\begin{smallmatrix}A_j\\~\\I_n\end{smallmatrix}\right)}_{Q_j}
{\upsilon}_j
=
\begin{smallmatrix}
\sum\limits_{j=1}^k{{\upsilon}}_j^*Q_j{\upsilon}_j
\end{smallmatrix}\succcurlyeq 0.
\end{matrix}
\]
\qed
\vskip 0.2cm

\noindent
In a similar way one can show a slightly stronger statement.

\begin{Cy}\label{Cy:MatrixConvex}
For given real parameters $\alpha$, $\beta$, $\gamma$ and $\delta$,
whenever not empty,
the family of all matrices $A\in\C^{n\times n}$ satisfying,
\begin{equation}\label{eq:BasicMatrixConvex}
\left(\begin{smallmatrix}A\\~\\I_n\end{smallmatrix}\right)^*
\left(\begin{smallmatrix}
{\alpha}I_n&&{\beta}I_n\\~\\{\gamma}I_n&&{\delta}I_n
\end{smallmatrix}\right)
\left(\begin{smallmatrix}A\\~\\ I_n\end{smallmatrix}\right)
=
{\begin{smallmatrix}
{\alpha}A^*A+{\beta}A^*+{\gamma}A+{\delta}I_n
\end{smallmatrix}}
\in\overline{\mathbf L}_{I_n}~,
\end{equation}
is matrix-convex.
\end{Cy}
\vskip 0.2cm

Indeed, by assumption,
\[
\left(\begin{smallmatrix}A\\~\\I_n\end{smallmatrix}\right)^*
\left(\begin{smallmatrix}
{\alpha}I_n&&{\beta}I_n\\~\\{\gamma}I_n&&{\delta}I_n
\end{smallmatrix}\right)
\left(\begin{smallmatrix}A\\~\\ I_n\end{smallmatrix}\right)
=
\begin{smallmatrix}
{\alpha}A^*A+{\beta}A^*+{\gamma}A+{\delta}I_n
\end{smallmatrix}
=
\underbrace{\scriptstyle{\alpha}A^*A+\frac{1}{2}
(\beta+\gamma)(A^*+A)+{\delta}I_n}_{\succcurlyeq 0}
+
\underbrace{\scriptstyle\frac{1}{2}(\beta-\gamma)(A^*-A)}_{\rm skew-Hermitian},
\]
so the claim follows from Lemma \ref{La:MatrixConvex}.
\vskip 0.2cm

\begin{center}
ACKNOWLEDGEMENT
\end{center}

The constructive reviews are well appreciated.

\end{document}